\documentclass[centertags,12pt]{amsart}
\usepackage{latexsym}
\usepackage{amssymb}

\numberwithin{equation}{section}
\makeatletter
\renewcommand{\subsection}{\@startsection
{subsection}{2}{0mm}{\baselineskip}{-0.25cm}
{\normalfont\normalsize\bf}}
\makeatother

\newtheorem{theorem}{Theorem}[section]
\newtheorem{proposition}[theorem]{Proposition}

\newtheorem{lemma}[theorem]{Lemma}

{\theoremstyle{remark}
\newtheorem{remark}[theorem]{Remark}}
\theoremstyle{definition}

\def\cC{\mathcal C}
\def\cD{\mathcal D}

\def\cX{\mathcal X}

\def\fq{\mathbf F_q}

\def\fr{\mathbf Fr}

\def\div{{\rm div}}

\def\deg{{\rm deg}}

\def\det{{\rm det}}

\begin{document}

\author[M.~Giulietti]{Massimo Giulietti}\thanks{2000 Math. Subj. Class.: Primary 14G,
Secondary 11G}
\thanks{This research was performed within the activity of GNSAGA of the Italian
CNR, with the financial support of the Italian Ministry MIUR, project
``Strutture geometriche, combinatorica e loro applicazioni", PRIN 2001-2002}

\address{Dipartimento di Matematica e Informatica \\ Universit\`a degli Studi di Perugia
\\ Via Vanvitelli, 1 \\ 06123 Perugia, Italy}

\title[On the number of rational points of a plane algebraic curve]{On
the number of ratonal points of a plane algebraic curve}

        \begin{abstract}
The number of $\fq$-rational points of a plane non-singular
algebraic curve $\cX$ defined over a finite field $\fq$ is
computed, provided that the generic point of $\cX$ is not an inflexion 
and that $\cX$ is Frobenius non-classical with respect to conics.
        \end{abstract}

\maketitle

{\bf Keywords:} Algebraic curves, Rational points, Frobenius non-classical curves.

\section{Introduction}
Let $\cX$ be a plane non-singular algebraic curve of degree $d$
defined over a finite field $\fq$ of order $q=p^r$ with $p$ prime.
Let $N$ denote the number of points of $\cX$ with coordinates in
$\fq$, also called $\fq$-rational points. The
Hasse-Weil theorem states that
$$ N=q+1-\sum_{i=1}^{2g} \alpha_i $$
where $\alpha_i$ are certain algebraic integers, and
$g=\frac{1}{2}d(d-1)$ is the genus of $\cX$. Nevertheless,
formulas for $N$ in terms of $d$ and some other projective
invariants of $\cX$ are only known for few curves, see \cite{gerard}, \cite{kristin}.
For instance, $N=q+1+(d-1)(d-2)\sqrt{q}$ for the Fermat curve $X^d+Y^d+1=0$ with
$\sqrt{q}\equiv -1 {\pmod d},$ and $q$ squared,  but this formula
does not hold true for $q=q_0^m$ with $m>2$ and
$q_0^{m-1}+\ldots+q_0+1\equiv -1 {\pmod d}$, see \cite{Ksz}.

In \cite{Ghefez-voloch} the authors pointed out that $N=d(q-d+2)$
when $\cX$ is Frobenius non-classical, that is the image $\fr(P)$
of a generic point $P$ of $\cX$ under the Frobenius map lies in
the tangent line at $P$. Note that for $p>2$ any Frobenius non-classical
curve is non-classical in the sense that every point of the curve
is an inflexion. An example of a Frobenius non-classical curve is
the Fermat curve of degree $d=\sqrt{q}+1$ with $q$ squared, also
called Hermitian curve.

In this paper we will be concerned with the case where
\begin{itemize}
\label{proper}
\item[A)] $\cX$ is classical;
\item[B)] $\cX$ is Frobenius non-classical with respect to conics,
that is $\fr(P)$ lies in the osculating conic $\cC_P$ to $\cX$ at
a generic point $P$ of $\cX$.
\end{itemize}

For $p\geq 5$, such a curve $\cX$ has the ``non--classical'' type
property that the intersection multiplicity $I(\cX,\cC_P;P)$  at a
generic point $P$ is a power $p^\nu$ of $p$. For $p^\nu=\sqrt{q}$,
an example of such a curve is the Fermat curve of degree
$\frac{1}{2}(\sqrt{q}+1)$ with $q$ squared. Our result is the
following theorem.
\begin{theorem}
\label{main} Let $\cX$ be a plane non-singular algebraic curve of
degree $d$ defined over a finite field of order $q=p^r$ with
$p\geq 5$ prime. Assume that $\cX$ satisfies conditions A) and B). If $d<p^\nu-1$, then
\begin{equation}
\label{formula}
N=\frac{1}{2}[d(q+5-2d)-k]
\end{equation}
where $k$ denotes the number of non-$\fq$-rational inflexion
points $P\in\cX$.
\end{theorem}

\section{Preliminary results}
An essential tool in the study of the number of $\fq$-rational
points of an algebraic curve defined over $\fq$ is the
St\"ohr-Voloch method. Here, we only summarize the results from
\cite[Sections 1-2]{Asv} which play a role in the proof of
Theorem \ref{main}.

Let $\cX$ be a plane non-singular algebraic curve of degree $d$
and genus $g$ defined over a finite field $\fq$ of order $q=p^r$,
with $p$ prime. Let $\bar{\mathbf F}_q(\cX)$ be the function field
of $\cX$, and $x_0,x_1,x_2$ $\fq$-rational functions in
$\bar{\mathbf F}_q(\cX)$ such that $\cX$ has homogeneous equation
$F(x_0,x_1,x_2)=0$. Also, assume that $\cX$ is classical, that is
the generic point of $\cX$ is not an inflexion.

The ramification divisor $R$ and the $\fq$-Frobenius divisor $S$ of $\cX$
are defined as follows
$$
R=\div (\det
\begin{pmatrix}
x_0 & x_1 & x_2 \\
D^{1}_t(x_0)& D^{1}_t(x_1) & D^{1}_t (x_2) \\
D^{2}_t(x_0)& D^{2}_t(x_1) & D^{2}_t (x_2)
\end{pmatrix}
)+3\div (dt)+3E\,,
$$
\begin{equation}\label{frob}
S=\div (\det \begin{pmatrix}
x_0^q & x_1^q &  x_2^q \\
x_0 & x_1 &  x_2 \\
D^{1}_t(x_0)& D^{1}_t(x_1) & D^{1}_t (x_2)
\end{pmatrix}
)+\div (dt)+(q+2)E\,,
\end{equation}
where $D_t^{(k)}$ is the $k$-th Hasse derivative with respect to a
separating variable $t$, and $E=\sum v_P(E)P$ with $v_P(E)=-\min
\{v_P(x_0),v_P(x_1), v_P(x_2)\}$. For $P\, \in \, \cX$, let $j(P)$
denote the intersection multiplicity of $\cX$ with its tangent
line at $P$.
\begin{proposition}\label{summary}
\begin{itemize}
\item[(a) ] $\deg R=3(2g-2)+3d$;
\item[(b) ] $\deg S=(2g-2)+(q+2)d$;
\item[(c) ] $v_P(R)\ge j(P)-2$; equality
holds if and only if $p$ does not divide
$j(P)(j(P)-1)/2$;
\item[(d) ] $v_P(S)\ge j(P)$
for $P$ $\fq$-rational point in $\cX$; equality
holds if and only if $p$ does not divide $j(P)-1$.
\end{itemize}
\end{proposition}

The order-sequence $(j_0(P),j_1(P),j_2(P),j_3(P),j_4(5),j_5(P))$
at $P\in\cX$ is defined to be the set of intersection
multiplicities of $\cX$ at $P$ with conics, arranged in increasing
order. Since $\cX$ is classical, two cases occur according as $P$
is an inflexion point or not, namely either
$(0,1,2,3,4,\epsilon(P))$ or $(0,1,2,j(P),j(P)+1,2j(P))$. Apart
from a finite number of points of $\cX$, we have the same order
sequence $(0,1,2,3,4,\epsilon)$. If $p\geq 5$ and condition B) is
satisfied, then $\epsilon=p^v$ for an integer $v \geq 1$, that is
$I(\cX,\cC_P;P)=\epsilon$ for the osculating conic $\cC_P$ at $P$.
To investigate the number of $\fq$-rational points on $\cX$ the
following result by Garcia and Voloch \cite{Egarcia-voloch} is
needed.
\begin{proposition}\label{Zprop5.2}
If $p\ge 5$ and $\cX$ if Frobenius non-classical with respect to
conics with order sequence $(0,1,2,3,4,p^\nu)$, then there exist
$\fq$-rational functions $z_0,z_1,z_2,z_3,z_4,z_5\in \bar{{\mathbf F}}_q(\cX)$ such that $$
z_0^{p^\nu}x_0^2+z_1^{p^\nu}x_0x_1+z_2^{p^\nu}x_0x_2+z_3^{p^\nu}x_1^2+
z_4^{p^\nu}x_1x_2+ z_5^{p^\nu}x_2^2=0\,. $$
\end{proposition}


\section{Curves which are Frobenius non-classical for conics}

Throughout this section we assume that $p\ge 5$  and that $\cX$
satisfies both conditions A) and B). Let $x:=x_1/x_0$,
$y:=x_2/x_0$ and $f(x,y)=0$ be a minimal equation of $\cX$. By
Proposition \ref{Zprop5.2}, there exist $\fq$-rational functions
$z_0,z_1,z_2,z_3,z_4,z_5\,\in\, {\bar {\mathbf F}}_q(\cX)$ such
that
\begin{equation}\label{Z1}
z_0^{p^\nu}+z_1^{p^\nu}x+z_2^{p^\nu}y+z_3^{p^\nu}x^2+z_4^{p^\nu}xy+
z_5^{p^\nu}y^2=0\,.
\end{equation}

Let $P=(a,b)$ be an affine point of $\cX$, and choose an index $j$
with $0\le j\le 5$ such that $v_P(z_j)\le v_P(z_i)$ for $0\le i
\le 5$. Putting $m_i=z_i/z_j\,\in \, {{\bar {\mathbf F}}_q}(\cX)$
we have $v_P(m_i)\ge 0$, and therefore $$
m_0^{p^\nu}+m_1^{p^\nu}x+m_2^{p^\nu}y+m_3^{p^\nu}x^2+
m_4^{p^\nu}xy+ m_5^{p^\nu}y^2=0\,, $$ with $m_j=1$. Let
$s_0=m_0(a,b)^{p^\nu}+m_1(a,b)^{p^\nu}x+m_2(a,b)^{p^\nu}y+m_3(a,b)^{p^\nu}x^2+
m_4(a,b)^{p^\nu}xy+ m_5(a,b)^{p^\nu}y^2\,\in\,{{\bar {\mathbf
F}}_q}(\cX)$. Then, $$
\begin{array}{rl}
s_0=& s_0-m_0^{p^\nu}-m_1^{p^\nu}x-m_2^{p^\nu}y-m_3^{p^\nu}x^2-
m_4^{p^\nu}xy- m_5^{p^\nu}y^2 \\
= & (m_0(a,b)-m_0)^{p^\nu}+(m_1(a,b)-m_1)^{p^\nu}x+(m_2(a,b)-m_2)^{p^\nu}y+
(m_3(a,b)-m_3)^{p^\nu}x^2\\
{} &+(m_4(a,b)-m_4)^{p^\nu}xy+(m_5(a,b)-m_5)^{p^\nu}y^2
\end{array}
$$ and hence $v_P(s_0)\ge \min_{0\le i\le
5}v_P((m_i(a,b)-m_i)^{p^\nu}) \ge p^\nu$. Moreover, as
$m_j(a,b)=1$, the equation
$m_0(a,b)^{p^\nu}+m_1(a,b)^{p^\nu}X+m_2(a,b)^{p^\nu}Y+m_3(a,b)^{p^\nu}X^2+
m_4(a,b)^{p^\nu}XY+ m_5(a,b)^{p^\nu}Y^2=0$ defines a conic. Then
the following result is obtained.
\begin{lemma}\label{Z2.1}
For an affine point $P=(a,b)$ of $\cX$, let $\cD_P$ the conic of
equation
$$m_0(a,b)^{p^\nu}+m_1(a,b)^{p^\nu}X+m_2(a,b)^{p^\nu}Y+m_3(a,b)^{p^\nu}X^2+
m_4(a,b)^{p^\nu}XY+ m_5(a,b)^{p^\nu}Y^2=0.$$ Then the intersection
multiplicity $I(\cX,\cD_P;P)$ is at least $p^\nu$.
\end{lemma}
\begin{proposition}\label{Zcor1} Let $P=(a,b)\,\in\, \cX$.
If $d< p^\nu-1$, then
\begin{itemize}
\item[\rm{i)} ]$\cD_P$ coincides with the osculating conic $\cC_P$ of $\cX$ at $P$;
\item[\rm{ii)} ]${\mathbf F}r(P)\,\in\,\cC_P$.
\end{itemize}
\end{proposition}
\begin{proof}

i) If $P$ is not an inflexion point, then the osculating conic
$\cC_P$ is the only conic having intersection multiplicity with
$\cX$ at $P$ more than $4$, and hence i) holds. If $P$ is an
inflexion point, then $j(P)\leq d$. Thus $j_4(P)\leq d+1$. Since
$d<p^\nu-1$, the osculating conic of $\cX$ at $P$ turns out to be
the only conic having intersection multiplicity with $\cX$ at $P$
at least $p^\nu$.

ii) First note that since $P$ is an arbitrarily chosen  point on $\cX$, condition B)
is not sufficient to prove the assertion. From \cite[Corollary 1.3]{Asv} an equation for $\cC_P$ is
\begin{equation}\label{Zconica}
{\rm det}
\begin{pmatrix}
1 & x & y & x^2 & xy & y^2 \\
1& x(P) & y(P) & x^2(P) & xy(P) &
y^2(P) \\
D^{j_1}_t1& D^{j_1}_tx(P) & D^{j_1}_ty(P) & D^{j_1}_tx^2(P) & D^{j_1}_txy(P) &
D^{j_1}_ty^2(P) \\
D^{j_2}_t1& D^{j_2}_tx(P) & D^{j_2}_ty(P) & D^{j_2}_tx^2(P) & D^{j_2}_txy(P) &
D^{j_2}_ty^2(P) \\
D^{j_3}_t1& D^{j_3}_tx(P) & D^{j_3}_ty(P) & D^{j_3}_tx^2(P) & D^{j_3}_txy(P) &
D^{j_3}_ty^2(P) \\
D^{j_4}_t1& D^{j_4}_tx(P) & D^{j_4}_ty(P) & D^{j_4}_tx^2(P) & D^{j_4}_txy(P) &
D^{j_4}_ty^2(P)
\end{pmatrix}
=0\, ,
\end{equation}
with $t$ local parameter at $P$ and $(j_0,j_1,j_2,j_3,j_4)=(j_0(P),j_1(P),j_2(P),j_3(P),j_4(P))$.
Then, since
$j_4\le \max\{4,j(P)+1\}\le d+1<p^\nu$,
from the minimality of the $\fq$-Frobenius orders (\cite[p. 9]{Asv})
the rational function
$$
{\rm det}
\begin{pmatrix}
1 & x^q & y^q & (x^q)^2 & x^qy^q & (y^q)^2 \\
1 & x & y & x^2 & xy & y^2 \\
D^{j_1}_t1& D^{j_1}_tx & D^{j_1}_ty & D^{j_1}_tx^2 & D^{j_1}_txy &
D^{j_1}_ty^2 \\
D^{j_2}_t1& D^{j_2}_tx & D^{j_2}_ty & D^{j_2}_tx^2 & D^{j_2}_txy &
D^{j_2}_ty^2 \\
D^{j_3}_t1& D^{j_3}_tx & D^{j_3}_ty & D^{j_3}_tx^2 & D^{j_3}_txy &
D^{j_3}_ty^2 \\
D^{j_4}_t1& D^{j_4}_tx & D^{j_4}_ty & D^{j_4}_tx^2 & D^{j_4}_txy &
D^{j_4}_ty^2
\end{pmatrix}
$$
is equal to $0$. Therefore $(a^q,b^q)$ satisfies equation (\ref{Zconica})
and the assertion follows.
\end{proof}

The following proposition is the key fact to prove Theorem \ref{main}.
\begin{proposition}\label{Zpropf}
If $d< p^\nu -1$,
then $j(P)\,\in\,\{2,(p^\nu+1)/2\}$
for any point $P\, \in \, \cX$.
\end{proposition}
\begin{proof}

Let $P=(a,b)$ be a point on $\cX$. Introduce a new affine frame
$(X^{\prime},Y^{\prime})$ such that $P$ is taken to the origin and
the tangent line of $\cX$ at $P$ to the $X^{\prime}$-axis. The
corresponding change of coordinate functions from $(x,y)$ to
$(\xi,\eta)$ is given by
\begin{equation}\label{Zchange}
\begin{array}{ll}
x = & m_{11}\xi+m_{12}\eta+a\,,\\ y = & m_{21}\xi+m_{22}\eta+b\,,
\end{array}
\end{equation}
for some $m_{ij}\in\bar{{\mathbf F}}_q$, $i,j=1,2$. Equation
(\ref{Z1}) is invariant under this transformation. To see this,
put, for $0\le i\le 5$, $$
\begin{array}{l}
z_i(x,y)=z_i(m_{11}\xi+m_{12}\eta+a,m_{21}\xi+m_{22}\eta+b)={\bar z_i}(\xi,\eta)\,,
\\
a=c^{p^\nu}, \quad b=d^{p^\nu}, \quad m_{ij}=n_{ij}^{p^\nu},
\,i,j=1,2,
\end{array}
$$
and write ${\bar z_i}={\bar z_i}(\xi,\eta)$. Then, with
$$
\begin{array}{l}
\zeta_0(\xi,\eta)={\bar z_0}+c{\bar z_1}+d{\bar z_2}+c^2{\bar z_3}+cd{\bar z_4}
+d^2{\bar z_5}\,,\\
\zeta_1(\xi,\eta)=n_{11}{\bar z_1}+n_{21}{\bar z_2}+2cn_{11}{\bar z_3}
+(cn_{11}+dn_{21}){\bar z_4}+2dn_{21}{\bar z_5}\,,\\
\zeta_2(\xi,\eta)=n_{12}{\bar z_1}+n_{22}{\bar z_2}+2cn_{12}{\bar z_3}
+(cn_{22}+dn_{12}){\bar z_4}+2dn_{22}{\bar z_5}\,,\\
\zeta_3(\xi,\eta)=n_{11}^2{\bar z_3}+n_{11}n_{21}{\bar z_4}+n_{21}^2{\bar z_5}\,,\\
\zeta_4(\xi,\eta)=2n_{11}n_{12}{\bar z_3}+(n_{12}n_{21}+n_{11}n_{22}){\bar z_4}+2n_{21}n_{22}{\bar z_5}\,,\\
\zeta_5(\xi,\eta)=n_{12}^2{\bar z_3}+n_{12}n_{22}{\bar z_4}+n_{22}^2{\bar z_5}\,,\\
\end{array}
$$ equation (\ref{Z1}) becomes
\begin{equation}\label{Z2}
\zeta_0^{p^\nu}+\zeta_1^{p^\nu}\xi+\zeta_2^{p^\nu}\eta+_3^{p^\nu}\xi^2+\zeta_4^{p^\nu}\xi\eta+
\zeta_5^{p^\nu}\eta^2=0.
\end{equation}

Since the tangent line to $\cX$ at $P^\prime=(0,0)$ has equation
$Y^{\prime}=0$, we have $v_{P^\prime}(\eta)=j(P)$. Let $v_{P^\prime}(\zeta_i)=k_i$,
$i=0,1,\ldots,5$. The left-hand side in (\ref{Z2}) is the sum of
six rational functions with valuations at $P^\prime$: $$
\begin{array}{cc}
v_{P'} (\zeta_0^{p^\nu})=k_0p^\nu,
\qquad &
v_{P'}(\zeta_1^{p^\nu}\xi)=k_1p^\nu+1,
\\
v_{P'}(\zeta_2^{p^\nu}\eta)=k_2p^\nu+j(P),
\qquad &
v_{P'}(\zeta_3^{p^\nu}\xi^2)=k_3p^\nu+2,
\\
v_{P'}(\zeta_4^{p^\nu}\xi\eta)=k_4p^\nu+1+j(P),
\qquad &
v_{P'}(\zeta_5^{p^\nu}\eta^2)=k_5p^\nu+2j(P).
\end{array}
$$ At least two of these values must be equal, and less than or
equal to the remaining four. Hence one of the following holds: $$
\begin{array}{lll}
(k_0-k_1)p^\nu=1,
\quad &
(k_1-k_3)p^\nu=1,
\quad &
(k_2-k_4)p^\nu=1,
\\
(k_0-k_3)p^\nu=2,
\quad &
(k_0-k_2)p^\nu=j(P),
\quad &
(k_2-k_5)p^\nu=j(P),
\\
(k_0-k_4)p^\nu=1+j(P),
\quad &
(k_0-k_5)p^\nu=2j(P),
\quad &
(k_1-k_2)p^\nu=j(P)-1,
\\
(k_3-k_4)p^\nu=j(P)-1,
\quad &
(k_1-k_4)p^\nu=j(P),
\quad &
(k_4-k_5)p^\nu=j(P)-1,
\\
(k_1-k_5)p^\nu=2j(P)-1,
\quad &
(k_2-k_3)p^\nu=2-j(P),
\quad &
(k_3-k_5)p^\nu=2j(P)-2,
\end{array}
$$
Since $d<p^\nu-1$ we have $1+j(P)<p^\nu$. This leaves just two
possibilities: $j(P)=2$, $j(P)=(p^\nu+1)/2$.
\end{proof}

\section{Proof of Theorem \ref{main}}

To estimate the number of $\fq$-rational points of $\cX$, we will
use a procedure similar to that in \cite{Ghefez-voloch}. To do
this, we go on to study the ramification divisor  $R$ and the
$\fq$-Frobenius divisor $S$ of $\cX$.
        \begin{lemma}\label{Zchiave}
If $d<p^\nu-1$, then
for a point $P\,\in \,\cX$
$$
v_P(R)=j(P)-2\, , \qquad v_P(S)=
\begin{cases}
j(P) & \text{if $P\,\in\, \cX(\fq)$},\\
0 & \text{if $P\, \notin \, \cX(\fq)$ and $j(P)=2$},\\
j(P)-1 & \text{if $P\, \notin \, \cX(\fq)$ and $j(P)>2$}.
\end{cases}
$$
        \end{lemma}
\begin{proof} From Proposition \ref{Zpropf} $j(P)=2$ or $j(P)=
(1/2)(p^\nu+1)$, hence
$v_P(R)=j(P)-2$ by (c) of Proposition \ref{summary}. Suppose now that
$P\, \in \, \cX(\fq)$. Since $p$ does not divide $(j(P)-1)$,
from (d) of Proposition \ref{summary} it follows
$v_P(S)=j(P)$. For $P\, \notin \, \cX(\fq)$, we distinguish two cases.

If $j(P)=2$, any degenerate conic meet $\cX$ at $P=(a,b)$
with multiplicity at most $4$, and therefore the osculating conic $\cC_P$ at
$P$ is irreducible.
Moreover, $(a^q,b^q)$ belongs to $\cC_P$ by ii) of Propoposition \ref{Zcor1}.
Then $v_P(S)=0$ since otherwise $(a^q,b^q)$ would belong to the tangent
line $l_P$ at $P$, and there would exist
too many intersections between $l_P$ and $\cC_P$.

Suppose now that $j(P)>2$. Note that the osculating conic $\cC_P$
is the tangent line $l_P$ counted twice. From Proposition
\ref{Zcor1} it follows that $(a^q,b^q)\,\in\, \cC_P$ and hence
$(a^q,b^q)\, \in \, l_P$. Now, from equation
(\ref{frob}), $v_P(S)=v_P((x-x^q)D_t^1y-(y-y^q)D_t^1x)=
v_P((x-x^q)dy/dt-(y-y^q)dx/dt)$, with a separating variable $t\,
\in \, \fq(x,y)$ such that $v_P(dt)=0$. Since $v_P(S)$ is not
invariant under all affine transformations but only for those
fixing the plane over $\fq$, it is necessary to see how $v_P(S)$
changes under an $\bar{\mathbf F}_q$-linear transformation. With
$x,y,\xi,\eta$ as in (\ref{Zchange}), $$
\begin{array}{ll}
(x-x^q)dy-(y-y^q)dx=&[(a-a^q)m_{21}-(b-b^q)m_{11}]d\xi
\\
{} &+[(a-a^q)m_{22}-(b-b^q)m_{12}]d\eta
\\
{} &+(m_{11}m_{22}-m_{12}m_{21})(\xi d\eta-\eta d\xi)
\\
{} &-(m_{11}\xi+m_{12}\eta)^q(m_{21}d\xi+m_{22}d\eta)
\\
{} &+(m_{21}\xi+m_{22}\eta)^q(m_{11}d\xi+m_{12}d\eta)\,.
\end{array}
$$ Now we let $\tau=t(\xi,\eta)\,\in\,{\bar {\mathbf F}}_q(\cX)$.
%
%
By letting $\xi'=d\xi/d\tau$ and $\eta'=d\eta/d\tau$, the following formula is arrived
at:
$$
\begin{array}{ll}
v_P(S)=& v_{P'}(\xi')+v_{P'}\{[(a-a^q)m_{21}-(b-b^q)m_{11}]\\
{}& +[(a-a^q)m_{22}-(b-b^q)m_{12}]\eta'/\xi'\\
{}& +(m_{11}m_{22}-m_{12}m_{21})(\xi\eta'/\xi'-\eta)\\
{}& - (m_{11}\xi+m_{12}\eta)^q(m_{21}+m_{22}\eta'/\xi')\\
{}& + (m_{21}\xi+m_{22}\eta)^q(m_{11}+m_{12}\eta'/\xi')\}\,.
\end{array}
$$
Note that $v_{P'}(\xi)=1$ and $v_{P'}(\eta)=j(P)$ are both prime to $p$ by Proposition
\ref{Zpropf}. Hence, $v_{P'}(\xi')=0$ and
 $v_{P'}(\eta')=j(P)-1$. Furthermore, $(a-a^q)m_{21}-(b-b^q)m_{11}=0$
and $(a-a^q)m_{22}-(b-b^q)m_{12}\neq 0$, as the line joining
$(a,b)$ and $(a^q,b^q)$ is the tangent line at $P$. Hence $v_P(S)=
v_{P'}(\eta'/\xi')=j(P)-1$.
\end{proof}

Now we can prove the main result of the paper.

{\bf Proof of Theorem \ref{main}}
\begin{proof}
The genus $g$ of $\cX$ is equal to
$\frac{(d-1)(d-2)}{2}$, hence $\deg(R)=3d(d-3)+3d$ and
$\deg(S)=d(d-3)+d(q+2)$. Therefore $d(q+5-2d)=\deg(S)-\deg(R)=
\sum_{P\in\cX}[v_P(S)-v_P(R)]$. Then the assertion
follows from Lemma \ref{Zchiave}.
\end{proof}
\begin{remark}\label{rem4.2}
Note that the proofs of Lemma \ref{Zchiave} and of Theorem \ref{main} depend on conditions A) and B), and
on the following two facts arising from $d<p^\nu-1$:
\begin{enumerate}
\item[1)] ${\mathbf F}r(P) \in \cC_P$ for every $P\in \cX$ (see Proposition \ref{Zcor1});
\item[2)] $p$ does not divide $j(P)(j(P)-1)$ for every $P\in\cX$ (see Proposition \ref{Zpropf}).
\end{enumerate}
Therefore, if $\cX$ fulfills the above two conditions together with A) and
B) then its number
of $\fq$-rational points is given by equation (\ref{formula}). This happens for
instance for the following Fermat curves $\alpha X^d+\beta Y^d=1$ defined over $\fq$ (see \cite[Thm. 2, Thm.
3]{Ggarcia-voloch1}):
\begin{enumerate}
\item[i) ] $d=(q-1)/2(p^r-1)$ with  $\alpha^2,\beta^2\in{\mathbf F}_{p^r}\setminus \{0\}$;
\item[ii) ] $d=2(q-1)/(p^r-1)$ with  $p\equiv 1 \pmod 4$ and 
$\alpha,\beta$ non-zero squares in ${\mathbf F}_{p^r}$.
\end{enumerate}
\end{remark}
\begin{remark}
An example of plane non-singular algebraic curve which satisfies conditions A) and B)
but not equation (\ref{formula}) is given by
the Fermat curve $\alpha X^d+\beta Y^d=1$, with $d=2(q-1)/(p^r-1)$ and $\alpha,\beta$ non-squares in ${\mathbf
F}_{p^r}$. The number $N$ of its $\fq$-rational points
is $\frac{1}{2}d(q-1+d-d\psi+2\psi)$, where $\psi=0$ for $r$ odd
and $p\equiv 3\pmod 4$, $\psi=1$ otherwise (see \cite[Example (viii)]{Ggarcia-voloch1}).
Equation (\ref{formula}) would instead give
$N=\frac{1}{2}d(q+2-2d+\psi)$. It is easily seen that such a curve does not
satisfy condition 1) in Remark \ref{rem4.2}.
\end{remark}

\section*{Acknowledgments}
The author would like to thank Prof. G. Korchm\'aros and Prof. F. Torres for their useful comments.

\end{document}